\documentclass[12pt]{article}

\usepackage{amsfonts, amsmath, amsthm, amssymb}

\newtheorem{thm}{Theorem}

\newtheorem{lem}{Lemma}

\newtheorem{prop}{Proposition}

\newtheorem{remarks}{Remarks}

\date{}

\newcommand{\ee}{\mathbb{E}}

\newcommand{\pp}{\mathbb{P}}

\def\<{\langle}
\def\>{\rangle}

\def\beq{\begin{equation}}
\def\deq{\end{equation}}

\def\bdef{\begin{defn}}
\def\ndef{\end{defn}}

\def\bthm{\begin{thm}}
\def\nthm{\end{thm}}

\def\bprop{\begin{prop}}
\def\nprop{\end{prop}}

\def\brmk{\begin{remarks}}
\def\nrmk{\end{remarks}}

\def\bexa{\begin{exa}}
\def\nexa{\end{exa}}

\def\blem{\begin{lem}}
\def\nlem{\end{lem}}

\def\bcor{\begin{cor}}
\def\ncor{\end{cor}}

\def\bexe{\begin{exe}}
\def\nexe{\end{exe}}

\def\bprf{\begin{proof}}
\def\nprf{\end{proof}}

\def\bdes{\begin{description}}
\def\ndes{\end{description}}

\def\dsp{\displaystyle}

\begin{document}
%
%
\vskip 3mm

\noindent THE CENTRAL LIMIT THEOREM FOR LS ESTIMATOR IN SIMPLE
LINEAR EV REGRESSION MODELS
\vskip 3mm

\vskip 5mm

\noindent Yu Miao $^1$, Guangyu Yang $^2$ and Luming Shen $^3$

\noindent $^1$ Department of Mathematics and Statistics

\noindent Wuhan University

\noindent Hubei, China \ 430072

\noindent College of Mathematics and Information Science

\noindent Henan Normal University

\noindent Henan, China \ 453007

\noindent yumiao728@yahoo.com.cn

\noindent $^2$ Department of Mathematics and Statistics

\noindent Wuhan University

\noindent Hubei, China \ 430072

\noindent study\_yang@yahoo.com.cn

\noindent $^3$ Science College

\noindent Hunan Agriculture University

\noindent Hunan, China\ 410028

\noindent shenluming\_20@163.com

\vskip 3mm
\noindent Key Words: LS estimator; EV regression models;
central limit theorem.
\vskip 3mm

\noindent 2000 MR Subject Classification: Primary 62F12, Secondary
60F05.

\vskip 3mm

\noindent ABSTRACT

In this paper, we obtain the central limit theorems for LS estimator
in simple linear errors-in-variables (EV) regression models under
some mild conditions. And we also show that those conditions are
necessary in some sense.
\vskip 5mm

\noindent 1.   INTRODUCTION

The problem for convergence control of estimator is important in
practice. The central limit theorem plays a main role for
researching large sample problems. There are numerous studies on the
convergence in distribution and various estimation of deviation
probability are known, see e.g. Ibragimov and Has'miniskii (1979),
Ibragimov and Radavicius (1981). Especially, for classical maximum
likelihood estimator, Ibragimov and Has'miniskii (1979) succeeded in
proving the central limit theorem of MLE under some regularity
conditions.

In many economic applications, individual observations are very
naturally categorized into mutually exclusive and exhaustive groups.
For example, individuals can be classified into cohorts and workers
are employees of a particular firm. The simplest grouping estimator
involves taking the means of all variables for each group and then
carrying out a group-level regression by OLS or weighted least
squares (if there are different numbers of observations in different
groups). This estimator has been called the efficient Wald estimator
(cf. Angrist (1991)).

It is well known that there is small sample bias in the efficient
Wald estimator because the sample group means do not equal the
underlying population means. Deaton (1985) proposed an
Errors-in-Variables (EV) estimator to correct the effect of sampling
errors.

The main purpose of this paper is to study the asymptotic normality
of LS estimator for an EV model. For simplicity of representation,
as Liu and Chen (2005), we restrict ourselves to the case of simple
linear model:
 \beq\label{10}
  \eta_i=\theta+\beta x_i+\varepsilon_i, \ \
\xi_i=x_i+\delta_i,\ \  1\leq i\leq n
 \deq
 with the following assumptions:

 \bdes

\item{(1)} $\theta$, $\beta$, $x_1$, $x_2$, $\cdots$ are unknown
constants (parameters);

\item{(2)} ($\varepsilon_1$, $\delta_1$), ($\varepsilon_2$,
$\delta_2$), $\cdots$ are independent identically distributed
(i.i.d.) random vectors, $\varepsilon_1$, $\varepsilon_2$, $\cdots$
are i.i.d., $\delta_1$, $\delta_2$, $\cdots$ are i.i.d., and
$$
\ee\delta_1=\ee\varepsilon_1=0,\ \ 0<Var(\delta_1)=\sigma_1^2,\ \
Var(\varepsilon_1)=\sigma_2^2<\infty;
$$

\item{(3)} $\xi_i$, $\eta_i$, $i=1, 2, \cdots$ are observable.

\ndes From (\ref{10}) we have
 \beq\label{11}
 \eta_i=\theta+\beta\xi_i+\nu_i,\ \
\nu_i=\varepsilon_i-\beta\delta_i,\ \  1\leq i\leq n.
 \deq

Recently the studies for EV model have attracted much attention:
see, for example, Gleser (1981) obtained some large sample results
of estimation in a multivariate "errors in variables" regression
model. Amemiya and Fuller (1984) discussed the estimation for the
multivariate errors-in-variables model with estimated error
covariance matrix. Cui (1997) proved the asymptotic normality of
M-estimates in the EV model and Liu and Chen (2005) gave the
consistency of LS estimator of linear EV regression model under
rather weakly conditions and discovered that while in the ordinary
linear regression (with the errors i.i.d., the error-variance is
finite and non-zero) the weak, strong and quadratic-mean consistency
of the regression estimates are equivalent, it is not so in the EV
case, i.e., the quadratic-mean consistency requires much more
stronger conditions.

Consider formally (\ref{11}) as a usual regression model of $\eta_i$
on $\xi_i$, we get the LS estimator of $\theta$ and $\beta$ as
\begin{equation}
\hat{\beta}_n=\frac{\sum_{i=1}^n(\xi_i-\bar{\xi}_n)
(\eta_i-\bar{\eta}_n)}{\sum_{i=1}^n(\xi_i-\bar{\xi}_n)^2},\ \ \ \
\hat{\theta}_n=\bar{\eta}_n-\hat{\beta}_n\bar{\xi}_n,
\end{equation}
where $\bar{\xi}_n=n^{-1}\sum_{i=1}^n\xi_i$, and other similar
notations, such as $\bar{\eta}_n$, $\bar{\delta}_n$, are defined in
the same way.

 Under some common assumptions, Liu and
Chen (2005) proved the consistency of LS estimators of simple linear
EV model, and obtained that both weak and strong consistency of the
estimator are equivalent, but it is not so for quadratic-mean
consistency. They also proved that the following condition is the
sufficient and necessary condition for $\hat{\beta}_n$ being strong
and weak consistent estimate of $\beta$:
$$\lim\limits_{n\rightarrow\infty}n^{-1}S_n=\infty,$$
where $S_n=\sum_{i=1}^n(x_i-\bar{x}_n)^2$.

\vskip 3mm

\noindent 2. THE CENTRAL LIMIT THEOREM FOR THE ESTIMATOR OF
PARAMETER $\beta$

The model to be studied can be exactly described as follows:
\begin{equation}\label{15}
{\it\bf (C):}
\begin{cases}
\eta_i=\theta+\beta x_i+\varepsilon_i, \xi_i=x_i+\delta_i, 1\leq i\leq n; \\
(\varepsilon_i, \delta_i), 1\leq i\leq n, i.i.d.;\\
  \ee\delta_1=\ee\varepsilon_1=0, \ee\delta_1^2=\sigma_1^2,
\ee\varepsilon^2_1=\sigma_2^2, 0<\sigma^2_1, \sigma^2_2<\infty.
\end{cases}
\end{equation}
Here $(\xi_i, \eta_i)$, $1\leq i\leq n$ are observable, while $x_i,
1\leq i\leq n$, $\theta$, $\beta$, $\sigma_1^2$, $\sigma_2^2$ are
unknown parameters. We are mainly interested in the convergence rate
of the estimators for the regression coefficient $\beta$ and
constant term $\theta$.

Before our works, we mention a key lemma.
\begin{lem}\label{lem20}(Feller (1971), p530) For each $n\geq 1$, let $X_{n,1}, \ldots, X_{n,n}$
be $n$ independent variables with distribution $F_{n,k}$, $1\leq
k\leq n$. Let $T_n=\sum_{k=1}^nX_{n,k}$. Suppose that
$\ee(X_{n,k})=0$, $\ee(T_n^2)=1$, and that
 \beq
 \sum_{k=1}^n\int_{|x|>r}x^2F_{n,k}(dx)\to 0
 \deq
 for each $r>0$. Then the distribution of $T_n$ tends to $N(0,1)$,
 which is a standard normal distribution.
\end{lem}
\begin{thm}\label{thm1}
Under model (\ref{15}), let $S_n$ denote $\sum_{i=1}^n(x_i-\bar{x}_n)^2$ and assume that the following
condition satisfied:
 \beq\label{18}
 \lim_{n\to\infty}\frac{n}{\sqrt{S_n}}=0.
 \deq
 In addition, let the following conditions
be fulfilled: there exists a constant $\alpha>0$, such that
$\ee|\varepsilon_i|^{2+\alpha}< \infty$, $\ee|\delta_i|^{2+\alpha}<
\infty$ and
 \beq\label{208}
\lim_{n\to\infty}\max_{1\leq i\leq n}
\frac{|x_i-\bar{x}_n|}{S_n^{1/2}}=0, \deq
  then $\hat{\beta}_n-\beta$ satisfies the
 asymptotic normality, i.e.,
 \beq
 \frac{\sqrt{S_n}}{\sqrt{Var(\varepsilon_1-\beta\delta_1)}}(\hat{\beta}_n-\beta)\xrightarrow{d} N(0,
 1),
 \deq
 where $N(0, 1)$ is a standard normal distribution.
\end{thm}

\begin{remarks}{\rm Since $S_{n}=\sum_{i=1}^{n}(x_{i}-\bar{x}_{n})^{2}$, and together with (6) and (7),
we have \beq n\ll\sqrt{S_{n}}\ll\max_{1\leq i\leq
n}{|x_{i}-\bar{x}_{n}|}^{2}\ll{S_{n}}, \deq where $"a_{n}\ll b_{n}"$
\mbox{denotes} $a_{n}\rightarrow\infty,
\,\,b_{n}\rightarrow\infty\,\,\mbox{and}\,\,\lim\limits_{n\to\infty}
\displaystyle\frac{a_{n}}{b_{n}}=0$. This implies that the LS
estimator possess good large-sample properties only when
$\{x_{n},n\geq1\}$ has greater dispersion. For instance, if
$\{x_{n},n\geq1\}$ is a sequence of i.i.d. random variables with
common distribution $N(0,1)$, the LS estimates have no consistency
and asymptotic normality. In Section 4, we will show that the
condition (6) is necessary in some sense.}
\end{remarks}

\begin{proof} By simple calculation, we have
\beq\label{16} \hat{\beta}_n=\beta +
\frac{1}{\sum_{i=1}^n(\xi_i-\bar{\xi}_n)^2}\left(\sum_{i=1}^n(\xi_i-\bar{\xi}_n)
\varepsilon_i-\beta\sum_{i=1}^n(x_i-\bar{x}_n)\delta_i
-\beta\sum_{i=1}^n(\delta_i-\bar{\delta}_n)^2\right). \deq
 For obtaining our result, we need the following lemmas.
\begin{lem}\label{lem21}
For any $r>0$, we have \beq\label{210}
 \lim_{n\to\infty}\pp\left(\frac{\sum_{i=1}^n(\delta_i-\bar{\delta}_n)^2}{\sqrt{S_n}}\geq
 r\right)=0.
 \deq
\end{lem}
\begin{proof}
For any $r>0$, because of the following relation,
$$\aligned
& \pp\Big(\frac{\sum_{i=1}^n\delta_i^2}{\sqrt{S_n}}\geq
 r\Big)
 \leq& \pp\Big(\frac{\sum_{i=1}^{ [\sqrt{S_n}]}(\delta_i^2-\sigma_1^2)-\sum_{i=n+1}^{[\sqrt{S_n}]}
 (\delta_i^2-\sigma_1^2)+n\sigma_1^2}{\sqrt{S_n}}\geq
 r\Big),
\endaligned$$
where $[S_n]$ denote the integer part of $S_n$ and
$\sigma_1^2=\ee\delta_1^2$, we have
$$\aligned
    & \pp\left(\frac{\sum_{i=1}^n(\delta_i-\bar{\delta}_n)^2}{\sqrt{S_n}}\geq
 r\right)\leq \pp\left(\frac{\sum_{i=1}^n\delta_i^2}{\sqrt{S_n}}\geq
 r\right)\\
 \leq& \pp\Big(\frac{\sum_{i=1}^{ [\sqrt{S_n}]}(\delta_i^2-\sigma_1^2)}{\sqrt{S_n}}\geq
   \frac{r}{2}\Big)
   +\pp\Big(\frac{-\sum_{i=n+1}^{[\sqrt{S_n}]}(\delta_i^2-\sigma_1^2)}{\sqrt{S_n}}\geq
    \frac{r}{2}-\frac{n\sigma_1^2}{\sqrt{S_n}}\Big)\\
 \triangleq & \Delta_1+\Delta_2.
\endaligned$$
Next we will control $\Delta_1$ and $\Delta_2$ respectively. By the
law of large numbers, $\Delta_1$ is convergent to zero in
probability. For enough large $n$, we have
$$\aligned
\Delta_2=&\pp\Big(\frac{-\sum_{i=n+1}^{[\sqrt{S_n}]}(\delta_i^2-\sigma_1^2)}{\sqrt{S_n}}\geq
    \frac{r}{2}-\frac{n\sigma_1^2}{\sqrt{S_n}}\Big)\\
\leq
&\pp\Big(\frac{-\sum_{i=n+1}^{[\sqrt{S_n}]}(\delta_i^2-\sigma_1^2)}{\sqrt{S_n}-n}\geq
    \frac{r}{4}\Big).
\endaligned$$
Since the same reasons with $\Delta_1$, $\Delta_2$ is convergent to
zero in probability.  From the discussion above, we obtain the
lemma.
\end{proof}

\begin{lem}\label{lem22}For any $r>0$, we have
 \beq\label{211}
 \lim_{n\to\infty}\pp\left(\left|\frac{\sum_{i=1}^n(\xi_i-\bar{\xi}_n)^2}{S_n}-1\right|\geq
 r\right)=0.
 \deq
\end{lem}
\begin{proof} In fact, we only need to prove the form $\dsp\pp\Big(\frac{\sum_{i=1}^n(\xi_i-\bar{\xi}_n)^2}{S_n}-1\geq
 r\Big)$, and the proof of the other form $\dsp\pp\Big(\frac{\sum_{i=1}^n(\xi_i-\bar{\xi}_n)^2}{S_n}-1\leq
 -r\Big)$ is similar.
 By simple calculation, we have
$$\aligned
\sum_{i=1}^n(\xi_i-\bar{\xi}_n)^2=
\sum_{i=1}^n(x_i-\bar{x}_n)^2+2\sum_{i=1}^n(x_i-\bar{x}_n)(\delta_i-\bar{\delta}_n)+\sum_{i=1}^n(\delta_i-\bar{\delta}_n)^2,
\endaligned$$
and
$$\sum_{i=1}^n(\xi_i-\bar{\xi}_n)^2-S_n\leq \frac{r}{2}S_n+\frac{2+r}{r}\sum_{i=1}^n(\delta_i-\bar{\delta}_n)^2.$$
Thereby, we get
$$\aligned
\pp\Big(\frac{\sum_{i=1}^n(\xi_i-\bar{\xi}_n)^2}{S_n}-1\geq
 r\Big)\leq& \pp\Big(\frac{(2+r)\sum_{i=1}^n(\delta_i-\bar{\delta}_n)^2}{rS_n}\geq
 \frac{r}{2}\Big).
\endaligned$$
Hence, Lemma 3 turns to be obvious due to Lemma \ref{lem21}
\end{proof}
\begin{lem}\label{lem23}
For any $r>0$, we have
 $$\lim_{n\to\infty}\pp\left(\left|\frac{\sum_{i=1}^n(\delta_i-\bar{\delta}_n)\varepsilon_i}{\sqrt{S_n}}\right|\geq
 r\right)=0.$$
\end{lem}
\begin{proof}For any $r>0$, we have
$$\pp\left(\left|\frac{\sum_{i=1}^n(\delta_i-\bar{\delta}_n)\varepsilon_i}{\sqrt{S_n}}\right|\geq
 r\right)\leq \pp\left(\frac{\sum_{i=1}^n(\delta_i-\bar{\delta}_n)^2}{\sqrt{S_n}}
 +\frac{\sum_{i=1}^n(\varepsilon_i-\bar{\varepsilon}_n)^2}{\sqrt{S_n}}\geq
 2r\right).$$
 So, by Lemma \ref{lem21}, our result is obtained.
\end{proof}
\begin{lem}\label{lem24}
We have
 $$\frac{\sqrt{S_n}}{\sqrt{Var(\varepsilon_1-\beta\delta_1)}}\frac{\sum_{i=1}^n(x_i-\bar{x}_n)(\varepsilon_i-\beta\delta_i)}{S_n}
 \xrightarrow{d}
 N(0,1).$$
\end{lem}
\begin{proof}For any $n\geq 1$ and $1\leq i\leq n$, denote
$$X_{n,i}:=\frac{(x_i-\bar{x}_n)}{\sqrt{S_nVar(\varepsilon_1-\beta\delta_1)}}(\varepsilon_i-\beta\delta_i).$$
Then for any $n\geq 1$, $\{X_{n,i}\}_{i=1}^n$ is a sequence of
independent random variables, and $\ee X_{n,i}=0$, $\sum_{i=1}^n\ee
X_{n,i}^2=1$. Next, we shall consider the following control: for any
$r>0$,
$$\aligned
&\sum_{i=1}^n\ee\big[|X_{n,i}|^2\textbf{1}_{\{|X_{n,i}|> r\}}\big]\\
\leq& \frac{1}{S_nVar(\varepsilon_1-\beta\delta_1)}\sum_{i=1}^n
\int_{\left\{\frac{|(x_i-\bar{x}_n)x|}{\sqrt{S_nVar(\varepsilon_1-\beta\delta_1)}}>
r\right\}}\frac{|(x_i-\bar{x}_n)x|^{2+\alpha}}{r^\alpha
[S_nVar(\varepsilon_1-\beta\delta_1)]^{\alpha/2}}dF(x),
\endaligned$$
where $F(x)$ denotes the distribution of
$\varepsilon_1-\beta\delta_1$. Therefore,
$$\aligned
&\sum_{i=1}^n\ee\big[|X_{n,i}|^2\textbf{1}_{\{|X_{n,i}|> r\}}\big]\\
\leq & \frac{\ee|\varepsilon_1-\beta\delta_1|^{2
+\alpha}}{r^{\alpha}(S_nVar(\varepsilon_1-\beta\delta_1))^{1+\alpha/2}}\sum_{i=1}^n
|x_i-\bar{x}_n|^{2+\alpha}\\
\leq & \frac{\ee|\varepsilon_1-\beta\delta_1|^{2
+\alpha}}{r^{\alpha}(Var(\varepsilon_1-\beta\delta_1))^{1+\alpha/2}}\max_{1\leq
i\leq n} \frac{|x_i-\bar{x}_n|^{\alpha}}{S_n^{\alpha/2}},
\endaligned$$
and by (\ref{208}), we obtain
$$\lim_{n\to\infty}\sum_{i=1}^n\ee\big[|X_{n,i}|^2\textbf{1}_{\{|X_{n,i}|> r\}}\big]=0.$$
From Lemma \ref{lem20}, we obtain our lemma.
\end{proof}

Now we give the proof of Theorem \ref{thm1}.
 By (\ref{16}), we have
 \beq\label{209}
\sqrt{S_n}(\hat{\beta}_n-\beta)=\sqrt{S_n}
\frac{\sum_{i=1}^n(\delta_i-\bar{\delta}_n)
\varepsilon_i+\sum_{i=1}^n(x_i-\bar{x}_n)(\varepsilon_i-\beta\delta_i)
-\beta\sum_{i=1}^n(\delta_i-\bar{\delta}_n)^2}{\sum_{i=1}^n(\xi_i-\bar{\xi}_n)^2}.
 \deq
 From Lemma \ref{lem22} and Lemma \ref{lem23}, it is easy to see
 that in probability,
$$\frac{S_n}{\sum_{i=1}^n(\xi_i-\bar{\xi}_n)^2}\frac{\sum_{i=1}^n(\delta_i-\bar{\delta}_n)
\varepsilon_i}{\sqrt{S_n}}\to 0,$$ and with the same reasons, in
probability
$$\frac{S_n}{\sum_{i=1}^n(\xi_i-\bar{\xi}_n)^2}\frac{\sum_{i=1}^n(\delta_i-\bar{\delta}_n)^2
}{\sqrt{S_n}}\to 0.$$ Because of Lemma \ref{lem22} and Lemma
\ref{lem24}, we have
$$\frac{S_n}{\sum_{i=1}^n(\xi_i-\bar{\xi}_n)^2}
\frac{\sum_{i=1}^n(x_i-\bar{x}_n)(\varepsilon_i-\beta\delta_i)}{\sqrt{S_n}}\xrightarrow{d}
N(0, Var(\varepsilon_1-\beta\delta_1)).$$
Therefore, the proof of
Theorem \ref{thm1} is completed.

\end{proof}

\vskip 3mm

\noindent 3.  THE CENTRAL LIMIT THEOREM FOR THE ESTIMATOR OF
PARAMETER $\theta$

In this section, we shall discuss the central limit theorem of
estimator of parameter $\theta$. LS estimate of $\theta$ is
 \beq
\label{20}
 \hat{\theta}_n=\bar{\eta}_n-\hat{\beta}_n\bar{\xi}_n.
 \deq
 Hence
 \beq
\label{21}
 \hat{\theta}_n-\theta=(\beta-\hat{\beta}_n)\bar{x}_n+
(\beta-\hat{\beta}_n)\bar{\delta}_n-\beta\bar{\delta}_n+\bar{\varepsilon}_n.
 \deq
  Liu and Chen (2005), have proved the following
  theorem, which gave the sufficient
and necessary condition for $\hat{\theta}_n$ being weak consistent
estimate of $\theta$.
\begin{thm} Under model (\ref{15}), a sufficient
and necessary condition for  the weak consistency of
$\hat{\theta}_n$ is
 \beq
 \label{22}
  \lim\limits_{n\rightarrow \infty}n\bar{x}_n(S^{*}_n)^{-1}=0,\ \
  S^{*}_n=\max(n, S_n).
 \deq
\end{thm}
\begin{thm}\label{thm2} Under model (\ref{15}), assume that the conditions of
the Theorem \ref{thm1} are satisfied. In addition, we assume
 \beq
\frac{S_n}{n\bar{x}_n^2}\to \infty,
 \deq
then $\hat{\theta}_n-\theta$
satisfies the
 asymptotic normality, i.e.,
  \beq
 \frac{\sqrt{n}}{\sqrt{Var(\varepsilon_1-\beta\delta_1)}}(\hat{\theta}_n-\theta)\xrightarrow{d}
 N(0,1).
  \deq
\end{thm}
\begin{remarks} {\rm Since $n\bar{x}_{n}^{2}\ll S_{n}$ and the definition of $S_{n}$, we have
$\dsp
\lim\limits_{n\to\infty}{\frac{\sum_{i=1}^{n}{x_{i}^2}}{S_{n}}=1}$
and (\ref{22}) satisfied. Furthermore, we notice that the
normalization of $\beta_{n}$ is related to $\sqrt{S_{n}}$, however,
the normalization of $\theta$ is related to $\sqrt{n}$. It is not
mysterious, since the parameter $\beta$ has coefficients $x_{i}$,
but $\theta$ not in EV models. And this also suggests that the proof
of Theorem \ref{thm2} might be easier than that of Theorem
\ref{thm1}. }
\end{remarks}
\begin{proof}It is easy to see,
$$\frac{\sqrt{n}}{\sqrt{Var(\varepsilon_1-\beta\delta_1)}}(\bar{\varepsilon}_n-\beta\bar{\delta}_n)
\xrightarrow{d} N(0,1).$$ Hence, it is enough to show in probability
$$\frac{\sqrt{n}}{\sqrt{Var(\varepsilon_1-\beta\delta_1)}}(\beta-\hat{\beta}_n)(\bar{x}_n+
\bar{\delta}_n)\to 0.$$ By Theorem \ref{thm1}, we only need to show
in probability,
$$\frac{\sqrt{n}}{\sqrt{S_n}}(\bar{x}_n+
\bar{\delta}_n)\to 0.$$  The law of large numbers, and together with
the condition (17), yields our results.
\end{proof}

\vskip 3mm

\noindent 4. FURTHER DISCUSSIONS

 Recalling that Liu and
Chen (2005) proved that the following condition is the sufficient
and necessary condition for $\hat{\beta}_n$ being strong and weak
consistent estimate of $\beta$:
$$\lim\limits_{n\rightarrow\infty}n^{-1}S_n=\infty.$$
However, in Theorem \ref{thm1}, we assume the condition (\ref{18}),
i.e.,
 \beq\label{40}
 \lim_{n\to\infty}\frac{n}{\sqrt{S_n}}=0.
 \deq
It seems that the condition (\ref{40}) is strong.
 In this section, we will analyze it and show that it is
necessary in some sense.

From the proof of Lemma \ref{lem21}, we know that the condition
(\ref{40}) is a technical condition to prove the following limit
 \beq\label{42}
\lim_{n\to\infty}\pp\left(\frac{1}{\sqrt{S_n}}\sum_{k=1}^n
\delta_k^2\ge r\right)=0.
 \deq

Before our discussion, we need mention the following result.

\begin{thm}\label{thm4} (Petrov, 1987) Let $\{X_k\}_{k\ge 1}$ is a sequence of independent random
variables with the distribution $V_k(x)$ and $\{a_n\}_{n\ge 1}$ is a
sequence of increasing positive real number. Then
$$
\frac{1}{a_n}\sum_{k=1}^nX_k\xrightarrow{\pp} 0,
$$
if and only if the following conditions are satisfied
 \bdes

 \item{$(i)$} $$\sum_{k=1}^n\int_{|x|\ge a_n}dV_k(x)\to 0;$$

 \item{$(ii)$} $$\frac{1}{a_n^2}\sum_{k=1}^n\left\{\int_{|x|< a_n}x^2 dV_k(x)-\left(\int_{|x|< a_n}x dV_k(x)\right)^2\right\}\to 0;$$

 \item{$(iii)$} $$\frac{1}{a_n}\sum_{k=1}^n\int_{|x|< a_n}x dV_k(x)\to 0.$$

\ndes

\end{thm}

Under our model, we can take $a_n=\sqrt{S_n}$, $X_k=\delta^2_k$,
$k\ge 1$ in Theorem \ref{thm4}. If
$$
 \frac{1}{\sqrt{S_n}}\sum_{k=1}^n \delta^2_k \xrightarrow{\pp} 0,
$$
then, by (iii) of Theorem \ref{thm4}, we have

$$
\frac{1}{\sqrt{S_n}}\sum_{k=1}^n\ee\delta^2_kI_{\{\delta_k^2<\sqrt{S_n}\}}=\frac{n}{\sqrt{S_n}}\ee
\delta^2_1I_{\{\delta^2_1<\sqrt{S_n}\}}\to 0.
$$
Since $0<\ee \delta^2_1<\infty$, we get
$$
\frac{n}{\sqrt{S_n}}\to 0.
$$
That is to say,  (6) is the necessary condition for Theorem 1.

\vskip 3mm

\noindent 5.  ACKNOWLEDGEMENTS

The authors wish to thank Prof. L.M. Wu of Universit\'e Blaise
Pascal and  Wuhan University, and Prof. D.H. Hu of Wuhan University
for their helpful discussions and suggestions during writing this
paper. The second author was partially supported by the National
Natural Science Foundation of China (Grant No: 10371092) and the
Foundation of Wuhan University. At last, the authors are very
grateful to the conscientious anonymous referee for his very serious
and valuable report. His suggestions and comments have largely
contributed to the Section 4.

\vskip 3mm

\noindent BIBLIOGRAPHY

\vskip 3mm

\noindent Akahira, M. and Takeuchi, K. (1981). {\it Asymptotic
efficiency of statistical estimators.} Springer-Verlag New York.

\vskip 3mm

\noindent Amemiya, Y. and Fuller, W.A. (1984). Estimation for the
multivariate errors-in-variables model with estimated error
covariance matrix. {\it Ann. Statist.}, {\bf 12(2)}, 497-509.

\vskip 3mm

\noindent Angrist, J.D. (1991). Grouped Data Estimation and Testing
in Simple Labor Supply Models. {\it J. Econometrics}, {\bf 47},
243-265.

\vskip 3mm

\noindent  Cui, H.J. (1997). Asymptotic normality of M-estimates in
the EV model.
 {\it J. Sys. Sci. and Math. Sci.}, {\bf 10(3)}, 225-236.
\vskip 3mm

\noindent Deaton, A. (1985). Panel data from a time series of
cross-sections. {\it J. Econometrics}, {\bf 30}, 109- 126.

\vskip 3mm

\noindent Feller, W. (1971). {\it An introduction to probability
theory and its applications. Second Edition.} John Wiley and Sons,
Inc. Vol 2.

\vskip 3mm

\noindent Fuller, W. (1987). {\it Measurement Error Models.} Wiley,
New York.

\vskip 3mm

\noindent Gleser, L. J. (1981). Estimation in a multivariate "errors
in variables" regression model: Large sample results. {\it Ann.
Statist.}, {\bf 9(1)}, 24-44.

\vskip 3mm

\noindent  Ibragimov, I.A. and Has'miniskii, R.Z. (1979). {\it
Statistical Estimation.} Springer-Verlag  New York.

\vskip 3mm

\noindent Ibragimov, I. and Radavicius, M. (1981). Probability of
large deviations for the maximum likelihood estimator. {\it Soviet
Math. Dokl.}, {\bf 23(2)}, 403-406.

\vskip 3mm

\noindent  Liu, J.X. and Chen, X.R. (2005). Consistency of LS
estimator in simple linear EV regression models. {\it Math. Acta.
Sci.}, {\bf 25B(1)}, 50-58.

\vskip 3mm

\noindent Petrov, V.V. (1987). {\it Limit Theorems for Sums of
Independent Random Variables}, Nauka, Moscow. (in Russian). English
translation, (1991). Oxford Univ. Press, Oxford. Chinese
translation, (1991). by Huang, K.M. and Su, C. Technology Press of
China, Hefei, China.

\vskip 3mm

\noindent Stuart, A., Ord, J. K. and Arnold, S. (1999). {\it
Kendall's Advanced Theory of Statistics. Vol.2A.} Arnold, London.

\vskip 1in

\vskip 4mm

\vskip 4mm


\end{document}